\documentclass[12pt]{amsart}
\usepackage[francais]{babel}
\usepackage{graphics}
\usepackage{graphicx}
\usepackage{pdfpages}

\title{Francis Bessi\`Ere \\ \small un regard original sur les fondements} 

\author{Labib Haddad}
\address{120 rue de Charonne, 75011 Paris, France}
\email{labib.haddad@wanadoo.fr}

\usepackage{amssymb}
\usepackage{amsmath}
\usepackage{endnotes}

\newcommand{\su}{\subsection*}
\newcommand{\head}{\section*}
\newcommand{\noi}{\noindent}

\newcommand{\Ž}{\'e}

\newcommand{\ˆ}{\`a}
\newcommand{\}{\`u}

\newcommand{\leqs}{\leqslant}

\newcommand{\geqs}{\geqslant}

\newcommand{\guil}{\guillemotleft}  
\newcommand {\guir}{\guillemotright }

\newcommand {\et}{\ \text{et}\ }

\newcommand {\si}{\ \text{si}\ }

\newcommand {\sinon}{\ \text{sinon}\ }

\newcommand{\lopar}{\noi \{$\looparrowright$ \ }

\begin{document}
\maketitle

\thispagestyle{empty}

\markboth{Haddad}{Bessi\re}

{\bf We  draw attention to a manuscript submitted to the HAL Open Archives by Francis Bessi\re, where he tries to base mathematics on a \emph{translative} theory that could be shown  consistant using  only finitist methods, thus bypassing the impossibility shown by G\šdel for \emph{deductive} theories, such as [ZF], for example.}

\

{\small \bf On voudrait attirer l'attention sur un manuscrit  d\Žpos\Ž aux archives ouvertes HAL par Francis Bessi\re o\ il essaie de fonder les math\Žmatiques \ˆ l'aide d'une th\Žorie \emph {traductive }dont on pourrait \Žtablir la coh\Žrence [autrement dit, la non contradiction] par des m\Žthodes finitistes, contournant ainsi
l'impossibilit\Ž d\Žmontr\Že par G\šdel pour les th\Žories \emph{d\Žductives}, telles que [ZF], par exemple.}

\

\centerline{$*$ \ $*$ \ $*$}

\

Francis Bessi\re (1931-2015) est un ancien \Žl\ve de l'\'Ecole normale sup\Žrieure de la rue d'Ulm, promotion 1950. Agr\Žg\Ž de physique en 1954, il entre au service des \'Etudes \Žconomiques g\Žn\Žrales d'EDF (\'Electricit\Ž de France) en 1957 puis, apr\s un passage au Centre pour la recherche \Žconomique et ses applications en 1969, il  revient finir sa carri\re \ˆ EDF, jusqu'en  1991. 

\

Je l'ai connu vers l'\‰ge de 19 ans, superficiellement\!;  je n'appr\Žciais pas ses choix  politiques du moment. Apr\s la sortie de l'\'Ecole, nous nous sommes perdus de vue. Il \Žtait physicien, moi, math\Žmaticien. Je ne me doutais gu\re qu'il s'\Žtait int\Žress\Ž de pr\s au probl\me des fondements. C'est seulement apr\s sa mort que j'ai appris l'existence du manuscrit (inachev\Ž) suivant qu'il avait laiss\Ž sur le site d'archives ouvertes HAL  :

\

{\bf Ne peut-on pas contourner l'impasse identifi\Že par G\šdel et prouver la consistance des math\Žmatiques usuelles ? 2013. $<$hal-00812912$>$.}

\

Le manuscrit comporte 75 pages. Voici les r\Žsum\Žs qu'il en fait, en fran\c cais et en anglais.

\su{R\Žsum\Ž} Le point de d\Žpart de ce travail est le fait que les c\Žl\bres th\Žor\mes de G\šdel concernent sp\Žcifiquement les \guil syst\mes hypoth\Žti\-co-d\Žductifs\guir \ - que je pr\Žf\re appeler des \guil th\Žories d\Žductives\guir. Il devrait ainsi \tre possible de contourner l'impasse rep\Žr\Že par G\šdel en consid\Žrant d'autres types de th\Žories (sous r\Žserve qu'elles ne soient pas \Žquivalentes \ˆ des th\Žories d\Žductives). Je d\Žfinirai donc des \guil th\Žo\-ries traductives\guir \ puis des \guil th\Žories asymptotiques\guir. Le chapitre 1 et les treize
\guil M\Žmos\guir \ qui le compl\tent en posent les bases et pr\Žsentent (M\Žmo 4) une th\Žorie M dont les Annexes A, B et C prouvent qu'elle est bien \Žquivalente aux \guil Maths usuelles\guir, i.e. \ˆ [ZF] plus l'axiome de fondation (AF) - mais sans qu'il soit utile de prendre en compte l'axiome du choix (AC). En s'appuyant sur la notion intuitive de \guil dictif\guir \ (M\Žmo 10), les chapitres 2 et 3 introduisent alors des concepts que je crois nouveaux et tout-\ˆ-fait intuitifs : les \guil superblocs\guir \ de lettres li\Žes, et les dictifs \guil \ utiles\guir \ et \guil sub-transitifs\guir. Le chapitre 4 s'attaque \ˆ une famille de th\Žories traductives
qui sont toutes consistantes et aboutissant \ˆ M. Le chapitre 5, non r\Ždig\Ž, esp\Žrait bien achever ce travail, mais mon \‰ge et d'autres facteurs incontr\™lables m'ont emp\ch\Ž d'y parvenir : \ˆ d'autres chercheurs de prendre le relais ... 

\

\su{Summary.} The departure point of this study is the fact that the well known G\šdel's theorems apply specifically to \lq\lq hypothetico-deductive systems" - which I will rather call \lq\lq deductive theories". Then it may be possible to turn around the dead-end identified by G\šdel, by trying to find other kinds of theories (with the condition that they shall not be equivalent to a deductive theory). So I will define \lq\lq translative theories" and \lq\lq asymptotic theories". Chapter 1 and the thirteen following \lq\lq Memos" lay
down their foundations and present (Memo 4) the theory M : the Annexes A, B and C will prove that this theory is equivalent to \lq\lq usual Mathematics", i.e. [ZF] plus the axiom of foundation (AF) - but with no necessity to take into account the axiom of choice (AC). Using always the intuitive notion of \lq\lq dictives" (Memo 10), the chapters 2 and 3 will introduce some new concepts: the \lq\lq super-blocks" of bounded letters, and the \lq\lq useful" and \lq\lq sub-transitive" dictives. The chapter 4 deals with
a family of translative theories which are all consistent and leading to M. The unwritten chapter 5 hoped to finish the present work, but my age and other uncontrollable factors prevented me from a real success\!: I hope that other researchers will take over ... (sic) 

\

\centerline{$*$ \ $*$ \ $*$}

\

Dans les lignes qui suivent, on s'efforce de donner un aper\c cu, le plus clair, et sobre, possible du travail de Bessi\re. Pour la clart\Ž, l'expos\Ž est fait d'une succession de petits paragraphes, assez courts. Pour la sobri\Žt\Ž, on a \Žlagu\Ž un grand nombre de d\Žfinitions, de notations, de d\Žtails, de remarques  et de commentaires qui ne semblaient pas essentiels. Cet aper\c cu n'est qu'un simple guide pour la lecture du texte lui-m\me, foisonnant et volumineux, du manuscrit.

\

\noi Deux difficult\Žs se pr\Žsentent : le tr\s grand nombre de d\Žfinitions introduites et les multiples renvois, crois\Žs, qui s'entrem\lent, entre 4 chapitres, 13 m\Žmos et 3 annexes. L'une des principales r\Žf\Žrences est, bien \Žvidemment, le petit livre pr\Žcieux de
Jean-Louis KRIVINE, \it Th\Žorie axiomatique des ensembles\rm, P.U.F, Paris, 1969.

\

\lopar Il faut prendre garde, \Žgalement, \ˆ ceci. Souvent, certains mots apparaissent dans le texte avant leurs d\Žfinitions lesquelles ne sont pas toujours faciles \ˆ rep\Žrer dans la suite du manuscrit (ainsi, par exemple, de {\it constante} ou {\it atomique}.\}

\

\noi Bessi\re cherche \ˆ \Žtablir par des m\Žtodes finitistes que les math\Žma\-tiques habituelles [il dit {\it usuelles}] ont une th\Žorie coh\Žrente [il dit {\it logiquement consistante}]. On ne peut pas le faire avec une th\Žorie d\Žductive classique (qui serait r\Žcursive et engloberait l'arithm\Žtique) car elle tomberait sous le couperet du th\Žor\me de G\šdel. En contraste, aux th\Žories d\Žductives, commun\Žment appel\Žes hypoth\Žtico-d\Žductives, il oppose les th\Žories traductives qu'il  introduit  comme on va essayer de l'expliquer. 

\

On commence par d\Žfinir des langages ${\sf A, B, C, D, E, F,}$ de plus en plus riches, et des th\Žories qui leur sont associ\Žes. La th\Žorie $\sf A$ est la plus squelettique. La th\Žorie $\sf B$ est la th\Žorie bool\Ženne. La th\Žorie $\sf C$ est la th\Žorie de la logique classique bivalente. La th\Žorie $\sf D$ est celle des dictifs [mot nouveau forg\Ž par Bessi\re].

\

\noi Tout au long de son texte, Bessi\re distingue soigneusement les entiers intuitifs des entiers formels. Les entiers intuitifs seront d\Žsign\Žs par des majuscules $H, K, N, \dots$. On s'en sert \ˆ l'\Žcole \Žl\Žmentaire lorsque l'on compte les \Žl\Žments d'une {\bf liste (intuitive) finie}. En regard, {\bf une suite} est index\Že par {\bf des entiers formels}.

\su{Langage logique} Un langage logique est un langage formel ayant $\Psi$ et $\implies$ parmi ses symboles et qui les utilise comme dans le langage $\sf E$ d\Žcrit ci-dessous.

\su{Transcription et traduction } Soient $\sf R$ et $\sf S$ des langages formels. Une {\bf transcription} de $\sf S$ dans $\sf R$ est une {\it proc\Ždure intuitive finie} qui, \ˆ chaque \Žnonc\Ž $E$ de $\sf S$, associe des  \Žnonc\Žs $E^R$ de $\sf R$, des transcriptions :ces \Žnonc\Žs $E^R$ sont,  en quelque sorte,  des codes de l'\Žnonc\Ž $E$  dans le langage $\sf R$. Une {\bf traduction} de $\sf S$ dans $\sf R$ est une transcription univoque : \ˆ chaque \Žnonc\Ž $E$ de $\sf S$, elle associe un unique \Žnonc\Ž $E^R$\ de $\sf R$\!; autrement dit, c'est une application de $\sf S$ dans $\sf R$,  au sens habituel.

\

\noi On dira qu'une traduction de 
$\sf S$ dans $\sf R$  est {\bf fid\le} lorsque $\sf R$  et $\sf S$ sont des langages logiques, que la traduction de $\Psi$ est $\Psi$ et la traduction de $E\implies F$ est $E^R \implies F^R$.

\su{Th\Žorie} Une th\Žorie est un langage logique o\ l'on a d\Žfini les seuls \Žnonc\Žs auxquels on s'int\Žresse et, parmi eux, les \Žnonc\Žs {\it clos}. \`A chaque \Žnonc\Ž $E$ est associ\Ž sa {\it n\Žgation} $\neg E$. Ce langage est accompagn\Ž d'une {\it m\Žthode de s\Žlection} d\Žfinie dans un m\Žtalangage convenable : les \Žnonc\Žs ainsi s\Žlectionn\Žs parmi les \Žnonc\Žs clos sont appel\Žs {\it th\Žor\mes} de la th\Žorie.

\su{Th\Žories traductives} Toute traduction d'un langage logique $\sf S$ dans une th\Žorie de r\Žf\Žrence $\sf R$ donn\Že d\Žfinit une th\Žorie traductive  dont la r\gle est la suivante : un \Žnonc\Ž T est un th\Žor\me de  $\sf S$  si et seulement si sa traduction est un th\Žor\me de $\sf R$.

\su{Th\Žorie $\sf A$} Le langage de la th\Žorie  $\sf A$ est r\Žduit aux seuls symboles $\Psi$ et $\Theta$.  Ses seuls \Žnonc\Žs sont $\Psi$ et $\Theta$, et son seul th\Žor\me est l'\Žnonc\Ž $\Theta$ dont la n\Žgation est $\Psi$.

\

\noi  [Bessi\re justifie son choix de  $\Psi$ et de $\Theta$ en faisant remarquer que $\Psi$ est l'initiale du grec antique $\psi\epsilon\nu\delta o\varsigma$ (lire : pseudos), qui signifie \guil faux, mensonger\guir, et $\Theta$ est l'initiale de $\theta\epsilon\omega\rho\eta\?\mu\alpha$, \guil th\Žor\me\guir \ bien s\žr.]

\su{Th\Žorie $\sf B$} Le langage de la th\Žorie bool\Ženne $\sf B$  est  le langage logique le plus simple, celui dont les seuls symboles sont $\Psi$ et $\implies,$ sans aucun \Žnonc\Ž atomique. Tous les \Žnonc\Žs de $\sf B$ sont d\Žclar\Žs clos. La r\gle d'inf\Žrence et trois  sch\Žmas d'axiomes de logique \Žl\Žmentaire d\Žfinissent sur $\sf B$ une th\Žorie d\Žductive et logique,  la  Th\Žorie bool\Ženne.

\su{L'inf\Žrence} Si $A$ et $A \implies C$ sont des th\Žor\mes, alors $C$ est un th\Žor\me.

\su{Les trois sch\Žmas d'axiomes} Ce sont

$\neg\neg F \implies F$,

$D \implies (E \implies D)$,

$(A \implies (B\implies C)) \implies ((A\implies B) \implies (A \implies C))$.

\su{Th\Žorie $\sf C$} Le langage $\sf C$ poss\de trois symboles :  $\Psi$, $\implies$, $\lambda$. Ses seuls \Žnonc\Žs atomiques, appel\Žs  {\bf variables}, sont produits par la r\gle suivante :  $\lambda\Psi$ est une variable et, si $v$ est une variable, alors $\lambda v$ est une variable. Les variables forment ainsi une liste (intuitive) infinie. Sur ce langage, on d\Žfinit la Th\Žorie logique $\sf C$ par la r\gle d'inf\Žrence et les trois sch\Žmas d'axiomes pr\Žc\Ždents. Les th\Žor\mes de cette th\Žorie sont, par d\Žfinition, les sch\Žmas logiques. On sait (depuis longtemps) que cette th\Žorie est coh\Žrente, compl\te et d\Žcidable.

\su{Th\Žorie $\sf D$} Le langage $\sf D$ poss\de 5 symboles : $\Psi,\implies,  \delta,\varsigma,\in.$ On y distingue les dictifs et les \Žnonc\Žs, d\Žfinis par les r\gles suivantes : le symbole $\varsigma$ est un dictif et, si (Dictif) est un dictif, alors  $\delta$(Dictif) est un dictif. De m\me, $\Psi$ est un \Žnonc\Ž ainsi que $\implies$(Enonc\Ž)(Enonc\Ž) et $\in$(Dictif)(Dictif). Pour chaque entier intuitif $K$, on d\Žsigne par $D_K$ le dictif  $\delta\delta\delta\dots\delta\delta\delta\varsigma$ form\Ž par un $\varsigma$   pr\Žc\Žd\Ž  $K$ fois par $\delta$. Les dictifs se pr\Žsentent ainsi comme une {\it liste index\Že par les entiers intuitifs}.

\

\noi On d\Žsigne encore par $\sf D$, la  th\Žorie traductive associ\Že \ˆ ce langage, ayant la th\Žorie $\sf B$ comme r\Žf\Žrence, munie de la traduction binaire suivante.  Soit $A$ un \Žnonc\Ž de la forme $\in D_K D_N$. On \Žcrit l'entier $N$ sous sa forme binaire, $a_0a_1a_2\dots$ [par exemple, l'entier 6 s'\Žcrit $011$, sous forme binaire]. On prend 
$$A^B = \Theta \si a_K = 1 \ , \ A^B = \Psi \sinon.$$
Cette traduction est fid\le, de sorte que la th\Žorie $\sf D$ est \Žgalement coh\Žrente, compl\te et d\Žcidable.

\su{Interpr\Žtation} On dit que le dictif $X$ est un \Žl\Žment du dictif $Z$ lorsque  $\epsilon X Z$ est un th\Žor\me. Tout dictif $Z$ s'interpr\te ainsi comme un ensemble intuitif fini de dictifs, le dictif $\varsigma = D_0$ \Žtant interpr\Žt\Ž comme l'ensemble vide. R\Žciproquement, tout ensemble intuitif fini de dictifs est l'interpr\Žtation du dictif $D_N$ o\ $N$ est la somme des $2^K$ pour lesquels $\in D_KD_N$ est un th\Žor\me. Cela donne un contenu pr\Žcis aux notions intuitives  de partie d'un dictif, de successeur, d'inclusison,  d'\Žgalit\Ž, de r\Žunion, d'intersection, etc, pour les dictifs. 

\

\lopar Disons que cette interpr\Žtion est calqu\Že sur un mod\le classique, bien connu, de la th\Žorie des ensembles finis\}

\

\noi En particulier,  on introduit une liste (intuitive) illimit\Že, utile dans les prochains d\Žveloppements, la liste des $P_N$, en prenant $P_0 = \varsigma$ puis, par r\Žcurrence, $P_{N+1} = $ ensemble des parties de $P_N$.

\su{Le langage $\sf E$} Son alphabet poss\de 6 symboles : $ \Psi,\implies,  \lambda, \varsigma,\in, \forall.$ C'est un langage pour les \Žnonc\Žs des math\Žmatiques habituelles. 

\

\noi  On y  d\Žfinit deux types de \guil textes\guir \ : les litt\mes et les \Žnonc\Žs.

\

\noi R\gles de formation.  Les seuls litt\mes et  \Žnonc\Žs sont ceux que l'on obtient par les deux r\gles de formations suivantes.

\

\noi (1) Les litt\mes.  Le texte $\lambda\varsigma$ est un litt\me et,  si $L$ est un litt\me, alors $\lambda L$ est un litt\me.

\

\noi La hauteur H d'un litt\me  est le nombre de ses $\lambda$, de sorte que l'on a toujours   $H\geqs 1$.  

\

\noi (2) Les \Žnonc\Žs. Le symbole $\Psi$ est un \Žnonc\Ž et, si $A$, $B$, $C$ sont des \Žnonc\Žs, et $L$, $M$ sont des litt\mes, alors , $\forall A$, $\in LM$, $\implies BC$ sont des \Žnonc\Žs.

\

\noi $\forall A$ est un \Žnonc\Ž universalis\Ž; sa port\Že est le sous-texte $A$.

\noi $\in LM$ est une appartenance; sa port\Že  est le sous-texte $LM$.

\noi $\implies BC$ est une implication; sa port\Že  est le sous-texte $BC$.

\

\noi Comme c'est l'usage, au lieu de $\in LM$, on peut \Žcrire $L\in M$. De m\me, au lieu de $\implies BC$, on peut \Žcrire $B\implies C$.

\su{Trascription du langage $\sf E$}

\noi Un {\bf saturateur} est une occurrence de $\Psi$ ou d'un litt\me. Tout \Žnonc\Ž finit donc par un saturateur, (un $\Psi$ ou un $\varsigma$).

\

\noi Un {\bf signe} est une occurrence de $\forall, \implies, \in$, ou d'un saturateur.

\

\noi Par d\Žfinition, la {\bf profondeur}  du premier signe d'un \Žnonc\Ž est nulle; elle augmente de 1 apr\s chaque $\forall$; elle diminue apr\s chaque saturateur $s$ du nombre de sous-\Žnonc\Žs universalis\Žs que ce saturateur $s$  compl\te.

\

\noi Dire qu'un signe d\Žpend d'un autre signe donn\Ž veut dire qu'il est dans sa port\Že.

\

\noi Un signe de profondeur $P$ d\Žpend d'une unique occurrence d'un \lq$\forall$' de profondeur $K$, pour chacun des $K$ tels que $0\leqs K < P$.  Au total, il d\Žpend ainsi de $P$ signes \lq$\forall$'.

\

\noi Le {\bf niveau}  $N$ d'un litt\me est \Žgal, par d\Žfinition, \ˆ sa hauteur $H$ moins sa profondeur $P.$ 

\

\noi Ainsi,
\guil lorsque N est positif, on dit que ce litt\me est une occurrence du {\bf litt\Žral positif} de niveau $N$. Lorsque N est n\Žgatif ou nul, on dit que c'est une occurrence du {\bf litt\Žral li\Ž} \ˆ l'unique \lq$\forall$' de profondeur $-N$ dont il d\Žpend. Ce Òlitt\Žral li\ŽÓ est donc identifi\Ž par la position de ce \lq$\forall$' dans l'\Žnonc\Ž consid\Žr\Ž. De plus, chaque signe \lq$\forall$' introduit un litt\Žral li\Ž (qui peut \tre : \lq\lq vide", ou encore : \lq\lq inoccup\Ž", i.e. sans aucune occurrence effective).\guir

\

\noi Comme  il est d'usage,  $(\exists x) E$   abr\ge $\neg( \forall x) \neg E$.

\

\noi Les lettres des math\Žmatiques habituelles transcrivent les litt\Žraux, pas les litt\mes.  On dit qu'une lettre est {\bf libre} lorsqu'elle transcrit un litt\Žral positif. On dit qu'elle est {\bf li\Že} lorsqu'elle transcrit un litt\Žral li\Ž. 

\

\noi Un \Žnonc\Ž clos est un \Žnonc\Ž dont toute lettre libre est une constante.

\su{Les objets math\Žmatiques} On dit que l'\Žnonc\Ž $D(x)$ d\Žfinit l'objet $x$ lorsque $x$ est sa seule lettre libre et que l'\Žnonc\Ž suivant est un th\Žor\me (d'existence et d'unicit\Ž) :
$$(\exists x) D(x) \et (\forall y) D(y) \implies y=x.$$
Lorsque l'on fait de $D(a)$ un axiome explicite, on appelle {\bf constante} toute lettre non li\Že qui figure dans un tel axiome.

\su{Pr\Žnexe et d\Žnexe} On sait qu'un \Žnonc\Ž peut se pr\Žsenter sous plusieurs formes \Žquivalentes. Parmi elles, il y a les formes {\bf pr\Žnexes} celles o\ tous les quantificateurs, $\forall$ et $\exists$, sont regroup\Žs \ˆ gauche, en d\Žbut d'\Žnonc\Ž. Ici, un {\bf \Žnonc\Ž maxi-d\Žnect\Ž} est, par d\Žfintion, un \Žnonc\Ž \Žcrit avec des quantificateurs $\forall$ dont chacun est repouss\Ž le plus \ˆ droite possible. Pour plus de d\Žtails, on se reportera au manuscrit de Bessi\re.

 \su{La th\Žorie d\Žductive $\sf M$}
 On d\Žsigne, habituellement, par [ZF] la th\Žorie d\Žductive classique de Zermelo-Fraenkel, par [ZFC] la th\Žorie [ZF] augment\Že de l'axiome du choix et par (AF) l'axiome de fondation [1, p. 50]. On ne peut pas \Žtablir la coh\Žrence de la th\Žorie [ZF] par des moyens finitistes : c'est l'un des r\Žsultats de G\šdel. On sait, cependant, que les trois th\Žories, [ZF], [ZFC]  et [ZF] + (AF) sont \Žquicoh\Žrentes, autrement dit, si l'une d'elles est coh\Žrente, les deux autres le sont aussi. 
 
 \
 
 \noi Par commodit\Ž, Bessi\re convient de s'en tenir \ˆ la th\Žorie [ZF] + (AF) comme repr\Žsentant des math\Žmatiques habituelles. D'autres choix sont possibles et pourraient tout aussi bien faire l'affaire. Il introduit, dans le langage  $\sf E$, une th\Žorie d\Žductive $\sf M$, ad hoc,  puis il \Žtablit, m\Žticuleusement,  qu'elle est \Žquivalente \ˆ la th\Žorie [ZF] + (AF) : tout th\Žor\me de l'une {\it est} un th\Žor\me de l'autre.  Pour faire court, nous omettrons tous les d\Žtails de la d\Žmonstration  en renvoyant au manuscrit lui-m\me, en particulier, aux trois Annexes A, B, C.
 
 \su{Le langage $\sf F$} Il a pour ambition de rassembler les deux langages $\sf D$ et $\sf E$ en un langage unique.  L'alphabet du langage $\sf F$ poss\de ainsi les 7 symboles $ \Psi,\implies,  \delta, \lambda,  \varsigma,\in, \forall.$
 Les deux langages  $\sf D$ et $\sf E$ sont des sous-langages de $\sf F$.  Par d\Žfinition, les {\bf substantifs} de $\sf F$ sont les litt\mes et les dictifs. On y distingue deux types d'\Žnonc\Žs {\bf atomiques} : les \Žnonc\Žs universalis\Žs  et les appartenances, respectivement de la forme $\forall$(\' Enonc\Ž) et $\in$(Substantif)(Substantif). On dit qu'un \Žnonc\Ž est {\bf clos} lorsqu'il est sans litt\Žraux positifs. 
 
 \
 
 \noi \`A chaque dictif non vide $Z$, on associe une traduction fid\le de $\sf F$ dans $\sf F$, comme suit. [Voir M\Žmo 12 B et chapitre 4.1]. Si  un \Žnonc\Ž atomique autonome $A$ est une appartenance, on prend pour $A^Z$ l'\Žnonc\Ž $A$ lui-m\me. Si $A$ est de la forme $(\forall x_1 x_2\dots x_H) R(x_1,x_2,\dots,x_H)$ o\ $R(x_1,x_2,\dots,x_H)$ ne commence plus par $\forall$, on prend pour $A^Z$ la conjonction des $N^K$ \Žnonc\Žs que l'on obtient en substituant, de toutes les mani\res possibles, \ˆ chaque $x_H, H = 1, 2,\dots K$, un \Žl\Žment de $Z$, ce que l'on pourra \Žcrire en bref :

 $A^Z \equiv {\bf ET}(R(T_1,T_2,\dots, T_K) \ | \ T_1,T_2,\dots, T_K \in Z)$ 
 
 o\ {\bf ET} d\Žsigne une conjonction (globale), comme pour le symbole 
 
 de r\Žunion $\bigcup$.
 
 \
 
 \noi On d\Žsigne par $\sf F^C$ le sous-langage des \Žnonc\Žs clos : il englobe $\sf D$ et les \Žnonc\Žs strictement clos de $\sf E$, autrement dit les \Žnonc\Žs de $\sf E$ o\ aucun litt\Žral positif ne figure. On d\Žsigne par $\sf F^D$ le sous-langage de $\sf F^C$ form\Ž des \Žnonc\Žs {\bf maxi-d\Žnect\Žs},

 \su{Th\Žorie asymptotique} Sur un langage $\sf L$, on se donne une liste illimit\Že, $\sf L_1$, $\sf L_2, \dots, \sf L_R, \dots$, de th\Žories ayant les m\mes \Žnonc\Žs clos. On d\Žsigne par $\sf L^*$ la th\Žorie (dite {\bf asymptotique}) d\Žfinie sur le langage $\sf L$ en d\Žcidant qu'un \Žnonc\Ž est un th\Žor\me de $\sf L^*$ si et seulement si c'est un th\Žor\me de chacune des th\Žories $\sf L_R$ \ˆ partir d'un certain rang donn\Ž $R$. Bien entendu, il suffit que chacune des th\Žories $\sf L_R$ de la liste soit coh\Žrente pour que la th\Žorie $\sf L^*$ le soit.

\;\;

\noi Bessi\re se pose la question suivante sans y r\Žpondre, laissant \ˆ chacun le soin de le faire (pour lui-m\me) : le recours \ˆ des th\Žories asymptotiques peut-il \tre encore consid\Žr\Ž comme \guil finitiste\guir ?

\su{Le chapitre 4} C'est le plus long des quatre chapitres r\Ždig\Žs. Il comporte 11 pages; il est dense et {\it tr\s technique}. Il est consacr\Ž \ˆ 
un type de traductions de $\sf F$ dans $\sf D$ qui sont donc des th\Žories traductives.  Apr\s avoir introduit un type de traductions fid\les de $\sf F^D$ dans $\sf D$, on les prolonge en traductions de $\sf F$ dans $\sf D$ puis on   \Žtudie les conditions qui permettent aux th\Žories traductives ainsi d\Žfinies sur $\sf F$ de v\Žrifier, chaque fois,  certains axiomes et sch\Žmas d'axiomes \guil math\Žmatiques\guir \ de la th\Žorie $\sf M$. Ce chapitre comporte 5 longs paragraphes.

\

Disons quelques mots de chacun d'eux.

\su{Un type de traductions fid\les de $\sf F^D$ dans $\sf D$} On se donne, une fois pour toutes, une liste de dictifs, $C_1, C_2,\dots$ que l'on appelle {\bf les conteneurs}. Par exemple, la liste des $P_N$ d\Žfinis ci-dessus en m\me temps que le langage $\sf D$, et on introduit une fonction intuitive $\mu$ qui d\Žfinit une liste strictement croissante d'entiers $\mu K$. On suppose que $\varsigma$ appartient \ˆ chacun des $C_K$ et que chaque $C_K$ est contenu dans $C_{\mu K}$.

\

\noi \`A Chaque \Žnonc\Ž $B$ de $\sf F^D$ de la forme $(\forall t_1t_2\dots t_N)R(t_1t_2\dots t_N)$ o\ $R$ n'est pas universalis\Ž et $N \geqs 1$, on choisit une entier intuitif positif $K$ et on traduit $B$ par 
\[B^{/K} \equiv{\bf ET}(R(T_1,T_2,\dots, T_K) \ | \ T_1,T_2,\dots, T_K \in C_K)\tag{1}\] 

comme on l'a fait ci-dessus pour le dictif $Z$, dans le paragraphe qui 

introduit le langage $\sf F$.

\

\noi Les traductions \Žtudi\Žes vont appliquer une formule analogue
\ˆ la formule (1) \ˆ chaque sous-\Žnonc\Ž universalis\Ž d'un \Žnonc\Ž de $\sf F^D$, mais en modifiant le choix de l'indice $K$ suivant le sous-\ŽŽnonc\Ž qu'on consid\re.

\

\noi \`A chaque indice $H$, on asscocie une {\bf \Žtape} de traduction de $\sf F^D$ dans lui-m\me que l'on d\Žsigne par $^{\bullet H}$ et qui laisse invariants, en particulier, les \Žnonc\Žs de $\sf D$. \`A chaque valeur de $K$, on associe la traduction de $\sf F^D$ dans $\sf D$, {\it compos\Že} des \Žtapes d'indices $K, \mu K, \mu\mu K, \dots$ jusqu'\ˆ obtenir un \Žnonc\Ž de $\sf D$, ce qui arrive n\Žcessairement, comme on le montre. Cette traduction est fid\le. On la d\Žsigne par $^{[K}$. Toute th\Žorie d\Žfinie sur $\sf F^D$ par une traduction de ce type est donc coh\Žrente, compl\te et d\Žcidable.

\su{Prolongement vers $\sf F^C$ et vers $\sf F$} On prend une th\Žorie $\sf R$ engendr\Že sur $\sf F^D$ par une traduction du type $^{[K}$ d\Žfinie ci-dessus. On la {\it prolonge} d'abord vers $\sf F^C$, en composant avec la traduction $^d$ de $\sf F^C$ vers $\sf F^D$ appel\Že {\bf maxi-d\Žnexion} que l'on d\Žfinit, au pr\Žalable [dans 2.3].  Cette traduction de $\sf F^C$ dans $\sf D$ est fid\le. On la d\Žsigne par $^{d[k}$. 

\

\noi On prolonge encore vers $\sf F$ en composant $^{d[K}$ avec une traduction $^u$ de $\sf F$ dans $\sf F^C$ [voir le d\Žtail dans le manuscrit]. On obtient une th\Žorie traductive que l'on d\Žsigne encore par $\sf F$. Bien que la traduction $^u$ ne soit pas fid\le, on montre que la th\Žorie  $\sf F$ est coh\Žrente, compl\te et d\Žcidable.

\

\noi On montre aussi que la r\gle d'inf\Žrence est valable dans cette nouvelle th\Žorie $\sf F$ et que les trois sch\Žmas d'axiomes de la th\Žorie $\sf B$ sont encore des sch\Žmas d'axiomes de $\sf F$.

\su{Le sch\Žma d'\Žgalit\Ž} Dans la th\Žorie $\sf M$, figure le sch\Žma d'axiome d'\Žgalit\Ž sous la forme suivante
\[(\forall a, b)((\forall x) (x\in a \iff x\in b)) \implies (P(a) \iff P(b))\tag{2}\]
On dit qu'un dictif  $X$  est {\bf sub-transitif} lorsque $A$ et $B$ \Žtant des \Žl\Žments de $X$, l'\Žgalit\Ž  $A\cap X = B\cap X$ implique $A=B$. 

\

\noi Il suffit de supposer que tous les conteneurs $C_K$ sont sub-transitifs pour que l'on ait la propri\Žt\Ž {\bf importante} suivante. Un \Žnonc\Ž du type (2) dont tous les dictifs appartiennent \ˆ $C_K$ est un th\Žor\me de la th\Žorie
traductive que $^{ud[K}$ d\Žfinit sur $\sf F$.

\su{Dictifs-tests, \Žl\Žments visibles et virtuels} Ce paragraphe est {\it techniquement tr\s laborieux}. Il est tr\s difficile d'en donner un aper\c cu. Il pr\Žpare les d\Žmonstrations du paragraphe suivant, le dernier.

\su{Un exemple : l'axiome de fondation} On y montre, en particulier ceci. Moyennant certaines conditions sur les conteneurs, l'axiome de fondation (AF) est un th\Žor\me de toute th\Žorie traductive sur $\sf F$ d\Žfinie  par une traduction $^{ud[K}$.

\

\centerline{$*$ \ $*$ \ $*$}

\

Bessi\re projetait de conclure son travail par un chapitre 5 qu'il a renonc\Ž \ˆ r\Ždiger apr\s la mort de son \Žpouse, Marie-Jeanne, \Žtant lui-m\me d\Žj\ˆ atteint d'une  d\Žg\Žn\Žrescence maculaire \ˆ l'\oe il droit.

\

\centerline {\bf L'essai de Bessi\re m\Žrite d'\tre repris, prolong\Ž et achev\Ž.}

\

\head{Bibliographie}

\

\noi [1] Jean-Louis KRIVINE, \it Th\Žorie axiomatique des ensembles\rm, P.U.F, Paris, 1969.

\

\

\

\

\

\

\enddocument